\documentstyle[11pt,mssymb]{article}
\def\up{\uparrow}
\def\Up{\Uparrow}
\def\upupup{\up\!\up\!\up}
\def\dn{\downarrow}
\def\Dn{\Downarrow}
\newlength{\squarewidth} 
\newlength{\diagwidth} 
\font\shmt=chess20
\def\sqskip{\hspace*\squarewidth}
\def\ltod#1{
  \if#1KJ\else
  \if#1kj\else
  \if#1QL\else
  \if#1ql\else
  \if#1RS\else
  \if#1rs\else
  \if#1BA\else
  \if#1ba\else
  \if#1NM\else
  \if#1nm\else
  \if#1SM\else
  \if#1sm\else
  \if#1PO\else
  \if#1po\else
  Z
  \fi\fi\fi\fi\fi\fi\fi\fi\fi\fi\fi\fi\fi\fi
  }
\def\llhack#1{                
  \if#10\sqskip\else          
  \if#1.\sqskip\else
  \if#1SN\else\if#1sn\else#1\fi\fi\fi\fi}
\def\evenrank#1#2#3#4#5#6#7#8{
  \llhack#1\ltod#2\llhack#3\ltod#4\llhack#5\ltod#6\llhack#7\ltod#8}
\def\oddrank#1#2#3#4#5#6#7#8{
  \ltod#1\llhack#2\ltod#3\llhack#4\ltod#5\llhack#6\ltod#7\llhack#8}
\def\board#1#2#3#4#5#6#7#8#9{
  \setlength{\unitlength}{1\squarewidth}
  \begin{picture}(9.5,9.8)
  {\shmt
  \put(.75,8.52){\evenrank #1}
  \put(.75,7.52){\oddrank  #2}
  \put(.75,6.52){\evenrank #3}
  \put(.75,5.52){\oddrank  #4}
  \put(.75,4.52){\evenrank #5}
  \put(.75,3.52){\oddrank  #6}
  \put(.75,2.52){\evenrank #7}
  \put(.75,1.52){\oddrank  #8}}
  \thinlines
  \put(.75,1.5){\framebox(8,8){}}
  \thicklines
  \put(.85,1.35){\line(1,0){8.07}}\put(8.9,9.4){\line(0,-1){8.07}}
  \put(.75,0.65){\makebox(8,1)[b]{#9}}
  \end{picture}
  }
\begin{document}

\vspace*{-6ex}

\centerline{\huge\bf On numbers and endgames:}

\centerline{\Large Combinatorial game theory in chess endgames}

\vspace*{4ex}

\centerline{\large Noam D. Elkies}

\vspace*{5ex}

{\small\em
Neither side appears to have any positional advantage
in the normal sense.  [\ldots]
the player with the move is able to arrange
the pawn-moves to his own advantage [and win] in each case.
It is difficult to say why this should be so,
although the option of moving a pawn one or two squares
at its first leap is a significant factor.}

--- Euwe and Hooper [EH, p.55], trying to explain why Diag.~5
(Example~87 in~[EH]) should be a first-player win

\vspace*{2ex}

{\em What shall we do with an Up?} 

--- Parker-Shaw [PS, 159--160 (Bass 2)]

\vspace*{2ex}

{\large\bf Introduction}

It was already noted in [WW, p.16] that combinatorial game theory
(CGT) does not apply directly to chess because the winner of a chess
game is in general not determined by who makes the last move, and indeed
a game may be neither won nor lost at all but drawn by infinite
play.\footnote{Of course infinite play does not occur in actual games.
Instead the game is drawn when it is apparent that neither side
will be able to checkmate against reasonable play.  Such a draw
is either agreed between both opponents or claimed by one of them
using the triple repetition or the 50-move rule.  These
mechanisms together approximate, albeit imperfectly, the principle
of draw by infinite play.}
Still, CGT has been effectively applied to other games such as
Dots-and-Boxes and Go which are not combinatorial games in the
sense of~[WW].  The main difficulty with doing the same for chess
is that the $8\times8$ chessboard is too small to decompose into
many independent subgames, or rather that some of the chess pieces
are so powerful and influence such a large fraction of the board's
area that even a decomposition into two weakly interacting
subgames (say a Kingside attack and a Queenside counteroffensive)
generally breaks down in a few moves.
Another problem is that CGT works best with ``cold'' games,
where having the move is a liability or at most an infinitesimal
boon, whereas the vast majority of chess positions are ``hot'':
Zugzwang\footnote{This word, literally meaning ``compulsion to move''
in German, has long been part of the international chess lexicon.}
positions (where one side loses or draws but would have done better
if allowed to pass the move) are already unusual, and positions of
mutual Zugzwang (henceforth mZZ) where neither side has a good
or even neutral move are much rarer.\footnote{Rare, that is,
in practical play; this together
with their paradoxical nature is precisely why Zugzwang and mZZ are
such popular themes in composed chess problems and endgame studies.}
This too is true of Go, but there the value of being on move, while
positive, can be nearly enough constant to be managed by
``chilling operators''~[GO], whereas the construction of
chess positions with similar behavior seems very difficult
and to date no such position is known.

To find interesting CGT aspects of chess we look to the endgame.
With most or all of the long-range pieces no longer on the board there
is enough room for a decomposition into several independent subgames.
Also it is easier to construct mZZ positions in the endgame because
there are fewer pieces that could make neutral moves; indeed it
is only in the endgame that mZZ ever occurs in actual play.
If a mZZ is sufficiently localized on the chessboard
one may add a configuration of opposing pawns which
must eventually block each other, and the first player who has no move
on that configuration loses the Zugzwang.  Furthermore this
configuration may split up as a sum of independent subgames.  The
possible values of these games are sufficiently varied that one may
construct positions that, though perhaps not otherwise intractable,
illustrate some of the surprising [ONAG] identities.  Occasionally
one can even use this theory to illuminate the analysis of a
chess endgame occurring in actual play.

We begin by evaluating simple pawn subgames on one file or two
adjacent files; this allows us to construct some novel mZZ positions
and explain the pawn endgame that baffled Euwe.  We then show
positions containing more exotic values: fractions, switches
and tinies, and loopy games.  We conclude with specific open
problems concerning the values which may be realized by positions
either on the $8\times8$ chessboard or on boards of other sizes.

{\bf Note:} In the vast majority of mZZ occurring in actual play
only a half point\footnote{In tournament chess a win, draw or loss
is worth 1, 1/2, or 0 points respectively.} is at stake:
one side to move draws, the other loses.  We chose to illustrate
this article with the more extreme kind of mZZ involving the full
point: whoever is to move loses.  This is mainly because it is
easier for the casual player to verify a win than a draw, though
as it happens the best example we found in actual play is also
a full-point mZZ.  The CGT part of our analysis applies equally
to similar endgames where only half a point hinges on the mZZ.


{\large\bf Simple subgames with simple values}

{\bf Integers.} Integer values, indicating an advantage in
spare tempo moves,\footnote{A ``tempo move'' (a.k.a.\ ``waiting move'')
is a move whose only effect is to give the opponent the turn.}
are easy to find.  An elementary example follows:

\vspace*{2ex}

\centerline{
\board
{........}
{p.p.p...}
{..P.P...}
{........}
{p...P...}
{P.P..Kp.}
{..P...Pk}
{........}
{\bf Diagram 1}
}

The Kingside is an instance of the mutual Zugzwang (mZZ)
known in the chess literature as the ``tr\'ebuchet'':
once either White or Black runs out of pawn moves he must move his
King, losing the g-pawn\footnote{We use ``algebraic notation'' for
chess moves and positions: the ranks of the chessboard are numbered
1 through 8 from bottom to top; the columns (``files'') labeled
a through h from left to right; and each square is labeled by
the file and rank it is on.  Thus in Diag.~1 the White King is on f3.
Pawns, which stay on the same file when not capturing, are named
by file alone when this can cause no confusion.} and the game.
Clearly White has one free
pawn move on the e-file, and Black has two on the a-file, provided
he does not rashly push his pawn two squares on the first move.
Finally the c-file provides White with four free moves (the maximum
on a single file), again provided Pc2 moves only one square at a
time.  Thus the value of Diagram~1 for White\footnote{Henceforth
all game and subgame values are from White's perspective, i.e.\
White is ``Left'', Black is ``Right''.} is $1-2+4=3$, and
White wins with at least two free moves to spare regardless of who
moves first.

\pagebreak

{\bf Infinitesimals.} 
Simple subgames can also have values that are not numbers,
as witness the b- and h-files in Diagram~2:

\vspace*{2ex}

\centerline{
\board
{........}
{........}
{.p..p...}
{.......p}
{.P......}
{....P...}
{.......P}
{........}
{\bf Diagram 2}
}

The b-file has value $\{0|0\}=*$; the same value (indeed an isomorphic
game tree) arises if the White pawn is placed on a3 instead of b4.
The e-file has value zero, since it is a mZZ; this
is the identity $\{*|*\}=0$.  The h-file on the other hand has
positive value: White's double-move option gives him the advantage
regardless of who has the move.  Indeed, since the only Black move
produces a $*$~position, while White to move may choose between~0
and~$*$, the h-file's value is $\{0,*|*\}=\,\up$.  [While $\up$ is
usually defined as $\{0|*\}$, White's extra option of moving to~$*$
gives him no further advantage; in [WW] parlance it is ``reversible''
(p.~64) and bypassing it gives White no new options.
This may be seen on the chessboard by noting that in
the position: White pawns h2,f4 vs.\ Black pawns h5,f7, Black to move
loses even if White is forbidden to play h3 until Black has played h4:
1\ldots h4 2.f5, or 1\ldots f6 2.f5.]

This accounts for all two-pawn positions with the pawns separated
by at most two squares on the same file, or at most three on adjacent
files.  Putting both pawnse on their initial squares of either
the same or adjacent files produces a mZZ (value zero).
This leaves only one two-pawn position to evaluate, represented by
the a-file of Diagram~3:

\pagebreak

\vspace*{-2ex}

\centerline{
\board
{........}
{..p.p...}
{p.......}
{........}
{..PpP...}
{.....p.k}
{P..P.P..}
{.....KB.}
{\bf Diagram 3}
}

From our analysis thus far we know that the a-file has value
$\{0,*\,|\!\up\}$.  Again we bypass the reversible $*$ option, and
simplify this value to $\Up\!\!*=\{0|\!\!\up\}$.  Equivalently,
Diagram~3 (in which the c-, d- and e-files are $\Dn\!\!*$, and the
Kingside is mZZ for chess reasons) is a mZZ, and remains so if White
is forbidden to play a2-a4 before Black has moved his a6-pawn.
This is easily verified:  WTM (White To Move) 1.c5 a5 and wins
by symmetry, or 1.a4 d3 2.a5(c5) e5 or 1.a3 d3 etc.;
BTM 1\ldots a5 2.c5 wins symmetrically, 1\ldots d3 2.c5
e5 (else 3.e5 and 4.a3) 3.c6 a5 4.a4, 1\ldots c6 2.e5 and 3.a3,
1\ldots c5 2.d3 e5 (e6 3.e5 a5 4.a4) 3.a3 a5 4.a4.

On a longer chessboard we could separate the pawns further.
Assuming that a pawn on such a chessboard still advances one square
at a time except for an initial option of a double move on its first
move, we evaluate such positions thus: a White pawn on a2 against
a Black one on a7 has value $\{0,*\,|\!\Up\!\!*\}=\,\upupup$;
against a Black Pa8, $\{0,*\,|\!\upupup\}=\,\up\!\up\!\up\!\up\!\!*$;
and by induction on $n$ a Black pawn on the $(n+4)$th rank yields
$n$ ups or $n$ ups and a star according as $n$ is odd or even,
provided the board is at least $n+6$ squares wide so the Black
pawn is not on its initial square.  With both pawns on their
initial squares the file has value zero unless the board has
width 5 or 6 when the value is $*$ or $*2$ respectively.  Of course
if neither pawn is on its starting square the value is 0 or $*$
depending on the parity of the distance between them, as in the b- and
e-files of Diag.~2.

\pagebreak

Our next Diagram illustrates another family of infinitesimally valued
positions.

\centerline{
\board
{........}
{......p.}
{...p...p}
{p.......}
{.p..p.PP}
{........}
{PP.PP...}
{........}
{\bf Diagram 4}
}

The analysis of such positions is complicated by the possibility
of pawn trades which involve entailing moves: an attacked pawn
must in general be immediately defended, and a pawn capture
parried at once with a recapture.  Still we can assign standard
CGT values to many positions, including all that we exhibit in
Diagrams in this article, in which each entailing line of play
is dominated by a non-entailing one (see again [WW, p.64]).

Consider first the Queenside position in Diagram~4.  White to move
can choose between 1.b3 (value~0) and 1.a4, which brings about $*$
whether or not Black interpolates the {\em en passant}\/ trade
1\ldots b:a3 2.b:a3.  White's remaining choice~1.a3 would
produce an inescapably entailing position, but since Black can
answer 1.a3 with b:a3 2.b:a3 this choice is dominated
by 1.a4 so we may safely ignore it.  Black's move a4
produces a mZZ, so we have $\{0,*|0\}=\;\up\!\!*$ [WW, p.68].
Our analysis ignored the Black move 1\ldots b3?, but 2.a:b3
then produces a position of value~1 (White has the tempo move 3.b4
a:b4 4.b3), so we may disregard this option since $1>\;\up\!\!*$.

We now know that in the central position of Diagram~4 Black need
only consider the move d5 which yields the Queenside position of value
$\up\!\!*$, since after 1\ldots e3?\ 2.d:e3 White has at least
a spare tempo.  WTM need only consider 1.d4, producing a mZZ whether
or not Black trades {\em en passant}, since 1.d3?\ gives Black the
same option and 1.e3??\ throws away a spare tempo.  Therefore
the center position is $\{0\,|\!\up\!\!*\}=\;\Up$ [WW, p.73].  Both
this and the Queenside position turn out to have the same value as
they would had the pawns on different files not interacted.  This
is no longer true if Black's rear pawn is on its starting square:
if in the center of Diag.~4 Pd6 is placed on d7 the resulting position
is mZZ (WTM 1.d4 e:d3 2.e:d3 d6; BTM 1\ldots d5 2.e3 or d6 2.d4),
not $*$.  Shifting the Diagram~4 Queenside up one or two squares
produces a position (such as the Diag.~4 Kingside) of value $\{0|0\}=*$:
either opponent may move to a mZZ (1.h5 or 1\ldots g6), and neither
can do any better, even with Black's double-move option: 1.g5 h5
is again~$*$, and 1\ldots g5 2.h:g5 h:g5 is equivalent to 1\ldots g6,
whereas 1\ldots h5?\ is even worse.  With the h4-pawn on h3, though,
the double-move option becomes crucial, giving Black an advantage
(value $\{*|0,*\}=\;\dn$, using the previous analysis to evaluate
1.h4 and 1\ldots g6 as~$*$).

\vspace*{5ex}

{\large\bf An example from actual play:}

{\bf Schweda-Sika, Brno 1929}

\vspace*{2ex}

We are now ready to tackle a nontrivial example from actual play:

\vspace*{1ex}

\centerline{
\board
{........}
{.p.....p}
{p.......}
{....k...}
{....Pp..}
{.....K.P}
{PP......}
{........}
{Diagram 5}
}

\vspace*{-1ex}

This position was the subject of our opening quote from [EH].
On the e- and f-files the Kings and two pawns are locked in
a vertical tr\'ebuchet; whoever is forced to move there first
will lose a pawn, which is known to be decisive in such an endgame.
Thus we can ignore the central chunk and regard the rest as a
last-mover-wins pawn game.

As noted in the Introduction,
the $8\times8$ chessboard is small enough that a competent
player can play such positions correctly even without knowing
the mathematical theory.  Indeed the White player correctly evaluated
this as a win when deciding earlier to play for this position, and
proceeded to demonstrate this win over the board.  Euwe and Hooper
[EH, p.56] also show that Black would win if he had the move from
the diagram, but they have a hard time explaining why such a position
should a first-player win --- this even though Euwe held both the world
chess championship (1935--7) and a doctorate in mathematics.\footnote{
See [HW].  Of course Euwe and Hooper did not have the benefit of CGT,
which had yet to be developed.}

Combinatorial game theory tells us what to do: decompose the
position into subgames, compute the value of each subgame,
and compare the sum of the values with zero.  The central chunk
has value zero, being a mutual Zugzwang (mZZ).  The h-file
we recognize as $\Dn\!\!*$.  The Queenside is more complicated, with
a game tree containing hot positions (1.a4 would produce $\{2|0\}$)
and entailing moves (such as after 1.a4 b5); but again it turns
out that these are all dominated, and we compute that
the Queenside simplifies to $\up$.  Thus the total value
of the position is $\up + \Dn\!\!* =\, \dn\!\!*$.  Since this is
confused with zero, the diagram is indeed a first-player win.
To identify the Queenside value as $\up$ we show that White's move
1.h4, converting the Kingside to~$\dn$, produces a mZZ, using
values obtained in the discussion of Diag.~4 to simplify the
Queenside computations.  For instance,
1.h4 a5 2.h5 a4 3.h6 b6 4.b4 wins, but with WTM again after 1.h4,
Black wins after 2.a4 a5, 2.a3 h5, 2.b4 h5 (mZZ) 3.a3/a4 b5/b6,
2.b3 a5 (mZZ+$\dn$), or 2.h5 a5 and 3.h6 a4 or 3.a4(b3) h6.
Black to move from Diag.~7 wins with
1\ldots a5, reaching mZZ after 2.h4 a4 3.h5 h6 or 2.a4(b3) h6,
even without using the \ldots h5 double-move option (since without it
the h-file is $*$ and $\up +\, * =\, \up\!\!*$ is still
confused with zero).

\vspace*{4ex}

{\large\bf More complicated values: fractions, switches and tinies}

\vspace*{1ex}

{\bf Fractions.} Fractional values are harder to come by; Diagram~6
shows two components with value $1/2$.  In the Queenside component the
c2 pawn is needed to assure that Black can never safely play b4; a pawn
on~d2 would serve the same purpose.  In the configuration d4,e4/d7,f6
it is essential that White's e5 forces a pawn trade, i.e.\ that
in the position resulting from 1.e5 f5?\ 2.d5 White wins, either
because the position after mutual promotions favors White or because
(as in Diag.~6) the f-pawn is blocked further down the board.  Each
of these components has the form $\{0,*|1\}$, but (as happened in
Diagrams~2 and~3) White's $*$~option gives him no further advantage,
and so each component's value simplifies to $\{0|1\}=1/2$.
Since the seven-piece tangle occupying the bottom right corner of
Diag.~6 not only blocks the f6-pawn but also constitutes a (rather
ostentatious) mZZ, and Black's h-pawn provides him a free move,
the entire Diagram is itself a mutual Zugzwang illustrating the
identity $1/2+1/2-1=0$.

\vspace*{1ex}

\centerline{
\board
{........}
{...p....}
{.p...p..}
{.......p}
{..pPP...}
{..P..k.b}
{P.P..pNP}
{.....K.Q}
{\bf Diagram 6}
\board
{........}
{..p....p}
{p....p..}
{..p.....}
{..p...PP}
{..P..p..}
{PP...Pk.}
{....KR..}
{\bf Diagram 7}
}

\vspace*{-1ex}

What do we make of Diagram~7 then?  Chess theory recognizes the
five-man configuration around f2 as a mutual Zugzwang (the critical
variation is WTM 1.Rh1 K:h1 2.Kd2 Kg1!\ 3.Ke3 Kg2 and Black wins
the tr\'ebuchet; this mZZ is akin to the Kingside mZZ of Diagram~3,
but there the d-pawns simplified the analysis).  Thus we need to
evaluate the three pure-pawn subgames, of which two are familiar:
the spare tempo-move of Pc7, and the equivalent of half a spare
tempo-move White gets from the upper Kingside.  To analyze
the Queenside position (excluding Pc7), we first consider that
position after Black's only move a5.  From that position Black
can only play a4 (value~1), while White can choose between a4 and a3
(values~0 and~$*$ respectively), but not 1.~b3?\ c:b3 2.a:b3 c4!\
and the a-pawn promotes.  Thus we find once more the value
$\{0,*|1\}=1/2$.  Returning to the Diagram~7 Queenside, we now
know the value~1/2 after Black's only move~a5.  White's moves
a3 and a4 produce 0 and~$*$, and 1.b3 can be ignored because the
reply c:b3 2.a:b3 c4 3.b4 shows that this is no better than 1.a3.
So we evaluate the Queenside of Diagram~7 as $\{0,*|1/2\}=1/4$,
our first quarter.  Thus the whole of Diagram~7 has the negative
value $1/2-1+1/4=-1/4$, indicating a Black win regardless of who
has the move, though with BTM the only play is 1\ldots a5!\ producing
a $1/2+1/2-1$ \ mZZ.

On a longer chessboard we could obtain yet smaller dyadic fractions
by moving the b-pawn of Diag.~6 or the Black a-pawn of Diag.~7
further back as long as this does not put the pawn on its initial
square.  Each step back halves the value.  These constructions
yield fractions as small as $1/2^{N-7}$ and $1/2^{N-6}$ respectively
on a board with columns of length~$N\geq8$.

{\bf Switches and tinies.}  We have seen some switches (games $\{m|n\}$
with $m>n$) already in our analysis of 4-pawn subgames on two files
such as occur in Diagrams~4 and~5.  We next illustrate a simpler family
of switches.

\vspace*{2ex}

\centerline{
\board
{........}
{p....p..}
{..p....p}
{P.p.....}
{P....P..}
{..P..P.P}
{..P....P}
{........}
{\bf Diagram 8}
\board
{......bK}
{.p...kP.}
{...p.pp.}
{P..p..P.}
{...p..P.}
{..PP....}
{..p.....}
{s.S.....}
{\bf Diagram 9}
}

\vspace*{-2ex}

In the a-file of Diagram 8 each side has only the move a6.
If Black plays a6 the pawns are blocked, while White gains
a tempo move with a6 (cf.\ the \hbox{e-file} of our first Diagram),
so the a-file has value $\{1|0\}$.  On the c-file whoever
plays c4 gets a tempo move, so that file gives $\{1|-1\}=\pm1$.
Adding a Black pawn on c7 would produce $\{1|-2\}$; in general
on a board with files of length~$N$\/ we could get temperatures
as high as $(N-5)/2$ by packing as many as $N-3$ pawns on a
single file in such a configuration.\footnote{For large enough $N$\/
it will be impossible to pack that many pawns on a file starting
from an initial position such as that of $8\times8$ chess, because
it takes at least $n^2/4 + O(n)$ captures to get $n$ pawns of the
same color on a single file.  At any rate one can attain temperatures
growing as some multiple of~$\sqrt N$.}

The f-file is somewhat more complicated: White's f5 produces
the switch $\{2|1\}$, while Black has a choice between f6 and f5
which yield $\{1|0\}$ and~0.  Bypassing the former option we find that
the f-file shows the three-stop game $\{2|1\|0\}$.  Likewise
$\{4|2\|0\}$ can be obtained by adding a White pawn on~f2, and
on a longer board $n+1$ pawns would produce $\{2n|n\|0\}$.
The h-file shows the same position shifted down one square,
with Black no longer able to reach 0 in one step.  That file
thus has value $\{2|1\|1|0\}$, which simplifies to the number 1
as may be seen either from the CGT formalism or by calculating
directly that the addition of a subgame of value $-1$ to the h-file
produces mZZ.

Building on this we may construct a few tinies and minies,
albeit in more contrived-looking positions than we have
seen thus far (though surely no less natural than the positions
used in~[FL]).  In the Queenside\footnote{Excluding the pawns
on a5 and b7 which I put there only to forestall a White defense
based on stalemate.} on Diagram~9 the Black pawn on~c2
and both Knights cannot or dare not move; they serve only to
block Black from promoting after \ldots d:c3.  That is Black's
only move, and it produces the switch $\{0|-1\}$ as in the
a-file of Diag.~8.  White's only move is 1.c:d4 (1.c4?\ d:c4
2.d:c4 d3 3.N:d3 Nb3, or even 2\ldots Nb3 3.N:b3 d3) which yields mZZ.
Thus the Queenside evaluates to $\{0\|0|-1\}$, i.e.\ tiny-one.
Adding a fourth Black d-pawn on d7 would produce tiny-two,
and on larger boards we could add more pawns to get tiny-$n$
for arbitrarily large $n$.
In the Kingside of Diagram~9 the same pawn-capture mechanism
relies on a different configuration of mutually paralyzing pieces,
including both Kings.  With a White pawn on g3 the Kingside
would thus be essentially the same as the Queenside with colors
reversed, with value miny-one; but since White lacks that pawn
the Kingside value is miny-zero, i.e.\ $\dn$ [WW, p.124].
Black therefore wins Diag.~9 regardless of whose turn it is
since Black's Kingside advantage outweighs White's Queenside edge.

\vspace*{5ex}

{\large\bf Some loopy chunks}

Since pawns only move in one direction any subgame in which only pawns
are mobile must terminate in a bounded number of moves.  Subgames
with other mobile pieces may be unbounded, or ``loopy'' in the
[WW] terminology (p.314); indeed unbounded  games must have closed
cycles (loops) of legal moves because there are only finitely many
distinct chess positions.

\pagebreak

\centerline{
\board
{........}
{.......p}
{..k....P}
{p...p...}
{P.K.P..P}
{........}
{........}
{........}
{\bf Diagram 10}
\board
{........}
{.......p}
{.......P}
{p.k.p...}
{P...P..P}
{..K.....}
{........}
{........}
{$\phantom.$\bf Diagram 10$'$}
}

Consider for instance the Queenside of Diagram~10,
where only the Kings may move.  Black has no reasonable options
since any move loses at once to Kb5 or Kd5; thus the Queenside's
value is at least zero.  White could play Kc3, but Black responds
Kc5 at once, producing Diag.~$10'$ with value $\leq0$ because
then Black penetrates decisively at b4 or d4 if the White King budges.
Thus Kc3 can never be a good move from Diag.~10.  White can also
play Kb3 or Kd3, though.  Black can then respond Kc5, forcing Kc3
producing Diagram~$10'$.  In fact Black might as well do this at once:
any other move lets White at least repeat the position with
2.Kc4 Kc6 3.Kb3(d3), and White has no reasonable moves at all from
Diag.~$10'$ so we need not worry about White moves after Kb3(d3).
We may thus regard 1.Kb3(d3) Kc5 2.Kc3 as a single move which is
White's only option from the Diag.~10 Queenside.  By the same argument
we regard the Diag.~$10'$ Queenside as a game where White has no moves
and Black has only the ``move'' 1\ldots Kb6(d6) 2.Kc4 Kc6, recovering
the Queenside of Diagram~10.  We thus see that these Queenside
positions are equivalent to the loopy games called {\bf tis} and
{\bf tisn} in [WW, p.322].\footnote{These translate to {\bf istoo}
and {\bf isnot} in American English.}  Since the Kingside has value~1,
Diag.~10 ($1+\bf tis$) is won for White, as is Diag.~$10'$
($1+{\bf tisn}=\bf tis$) with BTM, but with WTM Diag.~$10'$
is drawn after 1.h5 Kb6(d6) 2.Kc4 Kc6 3.Kb3(d3) Kc5 4.Kc3 etc.

\pagebreak

\centerline{
\board
{....k...}
{.....ppp}
{........}
{........}
{........}
{........}
{PPP.....}
{...K....}
{\bf Diagram 11}
\board
{........}
{........}
{.k......}
{.....ppp}
{PPP.....}
{......K.}
{........}
{........}
{\bf Diagram 12}
}

We draw our final examples from the Three Pawns Problem (Diagram~11).
See [HW] for the long history of this position, which was finally
solved by Sz\'en around 1836; Staunton devoted twelve pages of
his Handbook to its analysis [S, 487--500].\footnote{Thanks to
Jurg Nievergelt for bringing this Staunton reference to my attention.}
Each King battles the opposing three pawns.  Three pawns on adjacent
files can contain a King but (unless very far advanced)
not defeat it.  Eventually Zugzwang ensues, and one player must
either let the opposing pawns through or push his own pawns
when they can be captured.  As with our earlier analysis we allow
only moves which do not lose a pawn or unleash the opposing
pawns; thus the last player to make such a move wins.  The Three
Pawns Problem then in effect splits into two equal and opposite
subgames.  One might think that this must be a mZZ, but in fact
Diag.~11 is a first-player win.  Diagram~12 shows a crucial point
in the analysis, which again is a first-player win despite the
symmetry.  The reason is that each player has a check (White's
a5 or c5, Black's f4 or h4) which entails an immediate King move:
Black is not allowed to answer White's 1.a5+ with the Tweedledum
move [WW, p.4] 1\ldots h4+, and so must commit his King before
White must answer the pawn check.  This turns out to be sufficient
to make the difference between a win and a loss in Diagrams~11
and~12.

\pagebreak

\centerline{
\board
{.k......}
{.P......}
{..P.....}
{P.......}
{.....ppp}
{........}
{......K.}
{........}
{\bf Diagram 13}
}

Diagram 13 is a classic endgame study by J. Behting using this material
(\#61 in [1234], originally published in {\em Deutsche Schachzeitung}
1929).  After 1.Kg1!\ the Kingside shows an important mZZ:
BTM loses all three pawns after 1\ldots g3 2.Kg2 or
1\ldots f3/h3 2.Kf2/h2 h3/f3 3.Kg3, while WTM loses after 1.Kh2(f1) h3,
1.Kf2(h1) f3, or 1.Kg2 g3 when at least one Black pawn safely promotes
to a Queen.  On other King moves from Diag.~13 Black wins: 1.Kf2(f1)
h3 or 1.Kh2(h1) f3 followed by \ldots g3 and White can no longer hold
the pawns, e.g.\ 1.Kh1 f3 2.Kh2 g3+ 3.Kh3 f2 4.Kg2 h3+ 5.Kf1 h2.
Thus we may regard Kg1 as White's only Kingside option in Diag.~13.
Black can play either \ldots g3 reaching mZZ, or \ldots f3/h3+ entailing
Kf2/h2 and again mZZ but BTM; in effect Black can interpret the
Kingside as either 0 or $*$.  In the Queenside, White to move can only
play a6 reaching mZZ.  Black to move plays 1\ldots Ka7 or Kc7,
when 2.a6 Kb8 is mZZ; but the position after 1\ldots Ka7(c7)
is not itself a mZZ because White to move can improve on 2.a6
with the sacrifice 2.b8Q+!\ K:b8 3.a6, reaching mZZ with BTM.
The positions with the Black King on b8 or a7(c7) are then seen
to be equivalent: Black can move from one to the other and
White can move from either to mZZ.  Thus the Diag.~13 Queenside is
tantamount to the loopy game with infinitesmial but positive value
which is called ${\bf over}=1/\bf on$ in [WW, p.317].
White wins Behting's study with 1.Kg1!\ Ka7(c7) 2.b8Q+!\ K:b8 3.c6
reaching mZZ; all other alternatives (except 2.Kg2 Kb8 repeating the
positions) lose: 1.a6?\ g3!, or 2.a6?\ Kb8.  Since {\bf over} exceeds
$*$ as well as~0, White wins Diag.~13 even if Black moves first:
1\ldots Ka7 2.Kg1!\ Kb8 3.a6, 1\ldots g3 2.a6, or 1\ldots f3/h3+
2.Kf2/h2 Ka7 3.b8Q+ K:b8 4.a6 etc.

\pagebreak

\centerline{
\board
{.k......}
{........}
{........}
{P.P.....}
{.......p}
{.....p..}
{.P.....K}
{........}
{\bf Diagram 14}
}

The mZZ in the analysis of the Diag.~13 Queenside after 2.b8Q+ K:b8
3.a6 is the only mZZ involving a King and only two pawns.  In other
positions with a King in front of two pawns either on adjacent files
or separated by one file, the King may not be able to capture the
pawns, but will at least have an infinite supply of tempo moves.
Thus such a position will have value {\bf on} or {\bf off}
[WW, p.317 ff.]\ according as White or Black has the King.
For instance in the Kingside of Diagram~14 the White must not
capture on h4 because then the f-pawn promotes, but the King
can shuttle endlessly between h2 and h3 while Black may not
move (1.\ldots f2?\ 2.Kg2 h3+ 3.K:f2!\ and the h-pawn falls next).
If White didn't have the pawn on b2, the Queenside would likewise
provide Black infinitely many tempo moves and the entire Diagram would
be a draw with value ${\bf on} + {\bf off} = {\bf dud}$ [WW, p.318].
As it is White naturally wins Diag.~14 since Black will soon run
out of Queenside moves.  We can still ask for the value of the
Queenside; it turns out to give another realization of {\bf over}.
Indeed the Black King can only shuttle between b7 and b8 since
moving to the a- or c-file loses to c6 or a6 respectively,
and until the b-pawn reaches b4 White may not move his other pawns
since c6/a6 drops a pawn to Kc7/a7.  We know from Diag.~13 that
if the b-pawn were on b5 the Diag.~14 Queenside would be mZZ.
The same is true with that pawn on b4 and the Black King on b7:
WTM 1.b5 Kb8, BTM 1\ldots Kb8 2.b5 or 2.c6 Kc7 3.b5.  Thus
pawn on b4 and King on b8 give $*$, as do Pb3/Kb7, while Pb3/Kb8
is again mZZ.  From b2 the pawn can move to mZZ against either
Kb7 or Kb8 (moving to $*$ is always worse, as in the Diag.~13
Queenside), yielding a value of {\bf over} as claimed.  Positions
such as this one, which show an advantage of {\bf over} thanks
to the double-move option, are again known to chess theory;
see for instance endgame \#55 of~[1234] (H. Rinck, published
in {\em Deutsche Schachzeitung} 1913) which uses a different
pawn trio.  Usually, as in that Rinck endgame, the position is
designed so that White can only win by moving a pawn to the
fourth rank in two steps instead of one.

\vspace*{5ex}

{\large\bf Open problems}

\vspace*{2ex}

We have seen that pawn endgames can illustrate some of the
fundamental ideas of combinatorial game theory in the
familiar framework of chess.  How much of CGT can be found
in such endgames, either on the $8 \times 8$ or on larger boards?
Of course one could ask for each game value in [WW] whether it can be
shown on a chessboard.  But it appears more fruitful to
focus on attaining specific values in endgame positions.
I offer the following challenges:

\begin{description}
\item[Nimbers.]  Do $*2$, $*4$ and higher Nimbers occur
on the $8\times 8$ or larger boards?  We have seen already
that on a file of length~6 the position Pa2~vs.~Pa5 gives
$*2$, and Pa2~vs.~Pb6 does the same for files of length~7.
But these constructions extend neither to longer boards
nor to Nimbers beyond $*2$ and $*3 = * + *2$.
\item[Positive infinitesimals.]  We have seen how to construct
tiny-$x$ for integers $x\geq0$ (Diag.~9).  How about other $x$
such as 1/2 or $1\!\!\up$?  Also, do the higher ups
$\up^2,\up^3,\ldots$ [WW, p.277 and 321] occur?
\item[Fractions.] We can construct arbitrary dyadic fractions
on sufficiently large chessboard.  Does 1/8 exist on the
$8\times 8$ board?  Can positions with value 1/3 or other non-dyadic 
rationals arise in loopy chess positions?  (Note that thirds
do arise as mean values in~[GO] thanks to the Ko rule.)
\item[Chilled chess?] Is there a class of chess positions that
naturally yields to chilling operators as do the Go endgames
of~[GO]?
\end{description}

In other directions, one might also hope for a more systematic
CGT-style treatment of {\em en passant} captures and entailing
chess moves such as checks, captures entailing recapture,
and threats to capture; and ask for a class of positions
on an $N\times N$\/ board that bears on the computational
complexity of pawn endgames as [FL] does for unrestricted
$N\times N$\/ chess positions.

{\bf Acknowledgements.} This paper would never have been
written without the prodding, assistance and encouragement
of Elwyn Berlekamp.  In the Fall of 1992 he gave
an memorable expository talk on combinatorial game theory,
during which he mentioned that at the time CGT was not known
to have anything to do with chess.  I wrote a precursor of this
paper shortly thereafter in response to that implied (or at least
perceived) challenge.  The rest of the material was mostly developed
in preparation for or during the MSRI workshop on combinatorial
games, which was largely Berlekamp's creation.  I am also grateful
for his comments on the first draft of this paper, particularly
concerning [GO].

This paper was typeset in \LaTeX, using Piet Tutelaers' chess
fonts for the Diagrams.  The research was made possible in part
by funding from the National Science Foundation, the Packard
Foundation, and the Mathematical Sciences Research Institute.

\end{document}